\documentclass[a4paper, 12pt]{article}

\usepackage{amssymb, amsmath, bbm, amsthm, eufrak}

\usepackage{fullpage}

\usepackage{graphicx}

\usepackage{tikz}

\theoremstyle{plain}
\newtheorem{theorem}{Theorem}[section]
\newtheorem{lemma}[theorem]{Lemma}

\numberwithin{equation}{section} 

\theoremstyle{definition}

\theoremstyle{remark}

\def\R{{\mathbb R}}
\def\Z{{\mathbb Z}}

\def\Nnint{{\mathbb Z}_+}

\def\Ss{\mathcal{W}^{n,n+1}}

\def\Sx{\mathcal{W}^n}
\def\Sy{\mathcal{W}^{n+1}}

\def\Ssw{\mathcal{W}^{n,n+1}_0}
\def\Sswo{\mathcal{W}^{n,n}_0}
\def\Sxw{\mathcal{W}^n_0}
\def\Sxwk{\mathcal{W}_0}
\def\Syw{\mathcal{W}^{n+1}_0}

\def\Gen{\mathcal{A}}
\def\Genw{\mathcal{A}_0}
\def\Genwo{\mathcal{A}_0}

\def\XGT{\mathfrak{X}}
\def\XGTd{\mathfrak{X}}
\def\XGTw{\mathfrak{X}}

\def\Kone{\ensuremath{\mathbb{K}}}
\def\SympKone{\ensuremath{\mathbb{K}^0}}

\newcommand{\Prob}{\mathop{\mathbb{P}}\nolimits} 
\newcommand{\E}{\mathop{\mathbb{E}}\nolimits} 
\newcommand{\Indi}{\mathop{\mathbbm{1}}\nolimits}

\makeatletter
\def\Ddots{\mathinner{\mkern1mu\raise\p@
\vbox{\kern7\p@\hbox{.}}\mkern2mu
\raise4\p@\hbox{.}\mkern2mu\raise7\p@\hbox{.}\mkern1mu}}
\makeatother

 \author{Jon Warren  and Peter Windridge}

\title{Some Examples of Dynamics for Gelfand-Tsetlin Patterns}

\begin{document}

\begin{center}
 {\Large \bf Some Examples of Dynamics for Gelfand-Tsetlin Patterns}

\vspace{.1in}

\begin{tabular}{cc}
 \textbf{Jon Warren} & \textbf{Peter Windridge} \\
{\tt warren@stats.warwick.ac.uk} & {\tt p.windridge@warwick.ac.uk} \\
Department of Statistics, & Department of Statistics, \\
University of Warwick, & University of Warwick, \\
Coventry CV4 7AL, UK & Coventry CV4 7AL, UK \\
\end{tabular}
\end{center}

\begin{abstract}
We give three examples of stochastic processes in the Gelfand-Tsetlin 
cone in which each component evolves independently apart from a
blocking and pushing interaction.  These processes give rise to couplings 
between certain conditioned Markov processes, last passage times and 
exclusion processes.  In the first two examples, we deduce known identities 
in distribution between such processes whilst in the third example, the components of the process
cannot escape past a wall at the origin and we obtain a new 
relation.
\end{abstract}

{\bf Keywords: Gelfand-Tsetlin cone; conditioned Markov process; exclusion process; 
last passage percolation; random matrices} 

{\bf AMS 2000 Subject Classification:} Primary 60J25; Secondary: 60C05.

Submitted to EJP on May 2, 2009, final version accepted July 16, 2009.

\section{Introduction}
In \cite{MR1682248}, the authors Baik, Deift and Johansson show
that suitably rescaled, the law of the longest increasing subsequence of a 
uniformly chosen random permutation of $\{1,2,\ldots,n\}$ converges, as $n$ tends to infinity, to that of the \emph{Tracy-Widom distribution}.  The latter, first identified in 
\cite{tracywidom94}, describes the typical fluctuations of the 
largest eigenvalue of a large random Hermitian matrix from the 
\emph{Gaussian unitary ensemble} (see \cite{MR2166458} for a definition).
This somewhat surprising discovery has been followed by much research 
which has shown that the Tracy-Widom distribution also occurs as a limiting law in various other
models such as last passage percolation \cite{johansson2007mdm,MR1737991}, exclusion processes \cite{sasamoto07}, random tilings
\cite{MR1900323,MR1737991} and polynuclear growth \cite{MR2018275, praehoferspohn02}. See also the survey
\cite{johansson2005rma}.  

Eigenvalues of random matrices are closely related to multi-dimensional random walks 
whose components are conditioned not to collide.  In particular, both fall into a class of processes with determinantal correlation structure and exhibit pairwise repulsion at a distance.
On the other hand, models such as the exclusion process
are defined by local ``hard edged'' interactions rather 
than particles repelling each other remotely.
This paper is concerned with showing how it is possible to connect these two types of model
by coupling processes of one class with processes from the other.

In common with previous works in this area, we realise these couplings via 
the construction of a stochastic process in the \emph{Gelfand-Tsetlin} cone
\[
\Kone_n = \{ (x^1, x^2,\ldots,x^n) \in \R^1\times \R^2 \times \ldots \times \R^n : x^{k+1}_i \leq x_i^k \leq x^{k+1}_{i+1} \}.
\]

A configuration $(x^1, \ldots, x^n) \in \Kone_n$ is called a
\emph{Gelfand-Tsetlin pattern} and we may represent the 
\emph{interlacing} conditions $x^{k+1}_i \leq x_i^k \leq x^{k+1}_{i+1}$ 
diagrammatically as follows.

\begin{center}
\begin{tabular}{*{9}{c}}
 && 	& 	&$x_1^1$&	&	&& \\
 && 	&$x^2_1$&	&$x^2_2$&	&& \\
 &&$x^3_1$&	&$x_3^2$&	&$x_3^3$&& \\
&$\Ddots$&&	&$\vdots$&	&&$\ddots$& \\
$x^n_1$& $x^n_2$& $x^n_3$ & &$\ldots$&	&&$x^n_{n-1}$ &$x^n_n$ \\
\end{tabular}
\end{center}

Suitable processes in the Gelfand-Tsetlin cone appear naturally in several settings, 
for example the particle process associated with a random domino tiling of 
the Aztec diamond \cite{nordenstam-2008} and the eigenvalues of a GUE 
matrix and its minors \cite{MR1818248}.   In other cases, the 
process in $\Kone_n$ is not evident at first sight 
and must be constructed, for example see the recent studies of 
asymmetric simple exclusion processes \cite{MR2438811,MR2430639}.

Most frequently, dynamics for the process in $\Kone_n$ 
are constructed using a combinatorial procedure known as the Robinson-Schensted-Knuth 
algorithm (see O'Connell \cite{oconnell2003crw}).
With RSK dynamics, the $n(n+1)/2$ components of the process 
are driven by a noise with only $n$ degrees of freedom, leading to 
strong correlations between components.

In this paper we consider some alternative dynamics in which 
every component of the process evolves independently 
except for certain blocking and pushing interactions that ensures the process 
stays in $\Kone_n$.  This approach yields a new relation between 
an exclusion type process constrained by an impenetrable wall 
and a multi-dimensional random walk with components conditioned 
to neither become disordered nor jump over the wall.  
Dynamics of this type have previously been considered
by Warren \cite{warren2005dsb} for Brownian particles (see also Toth and Veto 
\cite{toth2007srb}), by Nordemstam in the context of shuffling domino 
tilings of the Aztec diamond \cite{nordenstam-2008}, by Borodin and 
Ferrari in the context of surface growth models \cite{borodin2008agr}.  
Analogous dynamics have also previously been studied in the context 
of growth models where they are known as 
Gates and Westcott dynamics, see Pr\"{a}hofer and Spohn
\cite{praehoferspohn02} for example.

\section{Description of dynamics and results}

From here on, we work exclusively with Gelfand-Tsetlin 
patterns with integer valued components and hence modify 
our definition of $\Kone_n$ to
\[
\Kone_n = \{ (x^1, x^2,\ldots,x^n) \in \Z^1\times \Z^2 \times \ldots \times \Z^n : x^{k+1}_i \leq x_i^k \leq x^{k+1}_{i+1} \}.
\]

\subsection{Poisson case}\label{s:poissonintro}

Our first example consists of a continuous time $\Kone_n$ valued Markov process 
$(\XGT(t); t \geq 0)$ that determines the positions of $n(n+1)/2$ interlaced particles on the integer 
lattice $\Z$ at time $t$.  The stochastic evolution of the pattern $\XGT$ 
is as follows.

Fix a vector of rates $q \in (0,\infty)^n$ and identify each particle 
with its corresponding component in $\XGT$.  The particle $\XGT_1^1$ 
 jumps rightwards at rate $q_1 > 0$, i.e. after an exponentially distributed
waiting time of mean $q_1^{-1}$.  The two particles, $\XGT^2_1, \XGT^2_2$ 
corresponding to the second row of the pattern each jump rightwards at rate $q_2$ independently of $\XGT^1_1$ and each other unless either
\begin{itemize}
 \item $\XGT_1^2(t) = \XGT_1^1(t)$, in which case any rightward jump of $\XGT^2_1$ is suppressed (\emph{blocked}), or
 \item $\XGT_2^2(t) = \XGT_1^1(t)$, in which case $\XGT_2^2$ will be forced to jump (\emph{pushed}) if $\XGT^1_1$ jumps.
\end{itemize}

In general, for $k > 1$ and $1 \leq j < k$, each particle $\XGT_j^k$ attempts to jump rightwards at rate 
$q_k$, and will succeed in doing so unless it is blocked by particle 
$\XGT^{k-1}_j$.  Particle $\XGT_k^k$ can always jump rightwards at
rate $q_k$ without impediment.  In addition, if $\XGT^{k-1}_{j} = \XGT^k_{j+1}$, particle $\XGT^k_{j+1}$ is pushed to the right when $\XGT^{k-1}_{j}$ jumps.  
This blocking and pushing ensures that $\XGT(t)$ remains in $\Kone_n$ for every $t \geq 0$.
We will show that for certain initial conditions on $\XGT(0)$,
the bottom layer of the pattern, $(\XGT^n(t); t \geq 0)$, is distributed as a multi-dimensional 
random walk with components conditioned not to become disordered (Theorem \ref{t:charlierresult}).

To describe the result more precisely, recall that for 
$z \in \Sx = \{z \in \Z^n : z_1 \leq z_2 \leq \ldots \leq z_n \}$, the \emph{Schur function} $S_z:\R^n \to \R$ 
can be defined (see for example \cite{MR1464693}) as a sum of 
geometrically weighted patterns,
\begin{equation}\label{e:schur}
S_z(q_1, \ldots, q_n) = \sum_{x \in \Kone_n(z)} w^q(x).
\end{equation}
The sum is over $\Kone_n(z) = \{x \in \Kone_n : x^n = z \}$, the set of all Gelfand-Tsetlin patterns $x = (x^1, \ldots, x^n) \in \Kone_n$ with bottom row $x^n$ equal to $z$ and the geometric weight function is
\[
 w^q(x) = \prod_{i=1}^n q_i^{|x^i| - |x^{i-1}|},
\]
where $|z| = \sum_{i=1}^d z_i$ for $z \in \R^d$ and $|x^0| = 0$.

This definition gives a natural probability mass function
on patterns $x \in \Kone_n(z)$,
\begin{equation}\label{e:distonKn}
 M_z(x) = \frac{w^q(x)}{S_z(q)}.
\end{equation}

Suppose that $(Z(t);\; t \geq 0)$ is an $n$-dimensional random walk 
in which component $i$ is, independently of the other components, a 
Poisson counting process of rate $q_i$.  
The function $h:\Sx \to \R$ defined by
\begin{equation}\label{e:pphn}
  h(x) = q_1^{-x_1} \ldots q_n^{-x_n} S_x(q).
\end{equation}
is harmonic for $Z$ killed at the first instant it leaves $\Sx$
(see \cite{konig2001ncr} for example).  Hence, $h$ may be used to define
a new process, $Z^\dagger$, with conservative $Q$-matrix on $\Sx$ defined by
\[
Q(x,x+e_i) = q_i \frac{h(x + e_i)}{h(x)} = 
\frac{S_{x + e_i}(q)}{S_{x}(q)},\quad 1 \leq i \leq n, \; x \in \Sx,
\]
where $e_i$ is the standard basis vector, and the other off diagonal rates 
in $Q$ are zero.

This Doob $h$-transform, $Z^\dagger$, may be interpretted as a version of 
$Z$ conditioned not to leave $\Sx$ and is closely related to the Charlier 
orthogonal polynomial ensemble (again see \cite{konig2001ncr}).

In section \ref{s:poisson} we prove the following result, obtained 
independently by Borodin and Ferrari by another method in 
\cite{borodin2008agr}.

\begin{theorem}\label{t:charlierresult}
If $(\XGT(t);\; t \geq 0)$ has initial distribution $M_z(\cdot)$
for some $z \in \Sx$ then $(\XGT^n(t);\; t \geq 0)$ is 
distributed as 
an $n$ dimensional Markov process with conservative $Q$-matrix
\[
 Q(x,x+e_i) = \frac{S_{x + e_i}(q)}{S_{x}(q)}\Indi_{[x+e_i \in \Sx]}, \quad 1 \leq i \leq n, \; x \in \Sx. 
\]
and all other off diagonal entries set to zero, started from $z$.
\end{theorem}

Note that from structure of the initial distribution and the construction of $\XGT$, 
this theorem implies that in fact \emph{every} row of the pattern is distributed 
as a conditioned Markov process of appropriate dimension and rates.

Theorem \ref{t:charlierresult} readily yields a coupling of the type discussed 
in the introduction -- 
the (shifted) left hand edge $(\XGT^1_1(t), \XGT^2_1(t) - 1, \ldots, \XGT^n_1(t)- n +1; t \geq 0)$ of $\XGT$ has the same ``hard edged'' interactions as an asymmetric exclusion process
(the particle with position $\XGT^k_1(t)- k +1$, $1 \leq k \leq n$ takes unit jumps rightwards at rate $q_k$ but is barred from occupying the same site as any particle to its right).  However, Theorem \ref{t:charlierresult} implies that $(\XGT^n_1(t); t \geq 0)$ has the same law as $(Z^\dagger_1(t); t \geq 0)$, the first component of the random walk $Z$ conditioned to stay in $\Sx$, when started from $Z^\dagger(0) = z$.
Further we observe that when $z = (0,\ldots,0)$, $M_z$ is concentrated on the origin and 
a version of the left hand edge can be constructed from the paths of $Z$ via 
$\XGT^{1}_1(t) = Z_1(t)$ and
\[
 \XGT^{k+1}_1(t) = Z_{k+1}(t) + \inf_{0 \leq s \leq t}\left(\XGT^{k}_1(s) - Z_{k+1}(s) \right), \quad 1 \leq k < n.
\]

Iterating this expression and appealing to Theorem \ref{t:charlierresult}, 
\begin{equation}\label{e:oconnellyorident}
 \left( Z^\dagger_1(t); t \geq 0 \right) \overset{\mathrm{dist}}{=} \left(\inf_{0 = t_0 \leq t_1 \leq  \ldots \leq t_n = t} \sum_{i=1}^n \left(Z_i(t_i) - Z_i(t_{i-1})\right); t \geq 0\right).
\end{equation}

This identity was previously derived by O'Connell and Yor in \cite{oconnell2001rnc} using a construction based on the Robinson-Schensted-Knuth correspondence.

\subsection{Geometric jumps}
For our second example we consider a discrete time process 
$(\XGTd(t);\;t \in \Nnint)$ (where $\Nnint$ is the set of non-negative integers) in $\Kone_n$ in which components make
independent geometrically distributed jumps perturbed by interactions
that maintain the interlacing constraints.

Let $q$ be a fixed vector in $(0,1)^n$ and update the pattern 
at time $t$ beginning with the top particle by setting $\XGTd^1_1(t+1) = \XGTd^1_1(t) + \xi$, where $\xi$ is a geometric random variable with mean $(1-q_1)/q_1$.  
That is, the top most particle always takes geometrically distributed jumps 
rightwards without experiencing pushing or blocking.

Suppose rows 1 through $k-1$ have been updated for some $k > 1$ and we wish 
to update the position of the particle corresponding to the $j^{th}$ component of the 
$k^{th}$ row in the pattern, $\XGTd_j^k$. 
If $\XGTd_{j-1}^{k-1}(t+1) > \XGTd_j^{k}(t)$, 
then $\XGTd_j^k(t)$ is pushed to an intermediate position $\tilde \XGTd_j^k(t) = \XGTd_{j-1}^{k-1}(t+1)$, while if $\XGTd_{j-1}^{k-1}(t+1) \leq \XGTd_j^{k}(t)$, 
no pushing occurs and $\tilde \XGTd_j^k(t) = \XGTd_j^k(t)$.  

\begin{figure}[h!]
\begin{center}
\label{f:discreteblockpush}

\begin{tikzpicture}
\draw (4,2) circle (0.1cm) node[above left] {$\XGTd_1^1(t_0)$};

\fill (8,2) circle (0.1cm) node[above right] {$\XGTd_1^1(t_0+1)$};
\draw[->, thick] (4.075,2.075) arc (120:60:3.85cm);

\draw (2,1) circle (0.1cm) node[below] {$\XGTd_1^2(t_0)$};
\fill (4,1) circle (0.1cm) node[below] {$\XGTd_1^2(t_0+1)$};
\draw[style=dotted] (4,1.9)-- (4,1);
\draw[->, thick] (2.075,1.075) arc (120:90:3.85cm) -- node[right]{block}(4,1.1);

\draw (6,1) circle (0.1cm) node[below] {$\XGTd_2^2(t_0)$};
\draw[style=help lines] (8,1) circle (0.1cm) node[below] {$\tilde \XGTd_2^2(t_0)$};
\draw[->,style=help lines] (6.1,1) --node[above]{push} (7.9,1);
\fill (10,1) circle (0.1cm) node[below] {$\XGTd_2^2(t_0+1)$};

\draw[->, thick] (8.075,1.075) arc (120:60:1.85cm);

\end{tikzpicture}

\end{center}
\caption{Example of blocking and pushing}
\end{figure}
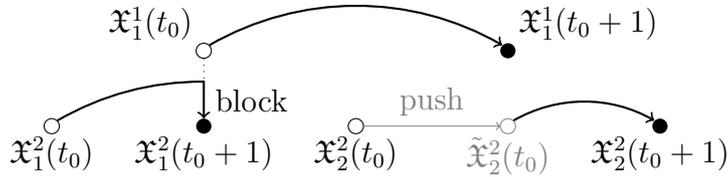

The particle $\XGTd_j^k$
then attempts to make a rightward jump of size that is geometrically distributed with mean $(1-q_k)/q_k$ from its intermediate position $\tilde  \XGTd_j^k(t)$ (so 
the particle is pushed \emph{before} it attempts to jump).  It always succeeds if $j = k$ (i.e. it is the right most particle) while if $j < k$, it cannot jump past $\XGTd_{j}^{k-1}(t)$, the position of particle to the right of it on the row above \emph{before} the update.
The leftmost particle $\XGTd_1^k$, $k > 1$ is not subject to pushing by any particle, but is still blocked by the ``ghost'' of the particle $\XGTd_{1}^{k-1}$.

To state the result, let us write $x \prec x^\prime$ when the inequality
$x_1 \leq x_1^\prime \leq x_2 \leq \ldots \leq x^\prime_{n-1} \leq x_n \leq x_n^\prime$ holds for $x,x^\prime\in \R^n$ and suppose $M_z$ is as defined in \eqref{e:distonKn}.  Then,

\begin{theorem}\label{t:discretegeometricintro}
If $\XGTd(0)$ has initial distribution $M_z(\cdot)$ for some $z \in \Sx$ then
$(\XGTd^n(t); t \in \Nnint)$ is distributed as an $n$ dimensional Markov process
in $\Sx$ with transition kernel
\[
p(x, x^\prime) = \prod_{i=1}^n (1-q_i) \frac{S_{x^\prime}(q)}{S_x(q)} \Indi_{[x \prec x^\prime]},\; x,x^\prime \in \Sx,
\]
beginning at $z$.
\end{theorem}

The Markov process with transition kernel $p$
can be described by a Doob $h$-transform - suppose $Z$ is now a discrete 
time random walk beginning at $z \in \Sx$
in which the $k^{th}$ component makes a
geometric($q_k$) rightward jump at each time step, independently of the other components.  Then the function $h$ defined in \eqref{e:pphn} is harmonic for $Z$ killed at the instant that the interlacing condition $Z(t) \prec Z(t+1)$ fails to hold (see \cite{oconnell2003crw}).  The corresponding $h$-transform $Z^\dagger$ is the discrete analogue of a process that arises from eigenvalues of Wishart matrices \cite{warrendiekerlargestevalue}.

The right hand edge of the pattern,  $(\XGTd^1_1(t), \XGTd^2_2(t), \ldots, \XGTd^n_n(t); t \in \Nnint)$ has a simple connection to the last passage percolation model with geometric weights that may be formulated as follows.  Suppose that $\eta_k(t)$ are independent geometric($q_k$) random variables attached to sites in the lattice $1 \leq k \leq n$, $t \geq 1$.  An increasing path $\pi$ from $(1,1)$ to $(t,k)$ is a collection of sites $\{(t_1, k_1), \ldots, (t_N, k_N)\}$, $N = t + k - 2$, such that the step $(t_{m+1}, k_{m+1}) - (t_m, k_m) \in \{(1,0),(0,1)\}$, and we denote the set of such paths by $\Pi(t,k)$.  The quantity of
interest is the $k$-dimensional process of last passage times
\[
 G_k(t) = \max_{\pi \in \Pi(t,k)} \sum_{(i,j) \in \pi} \eta_j(i),\; t \in \Nnint.
\]
It is not difficult to confirm that $(G_1(t), \ldots, G_n(t); t \in \Nnint)$ has
the same law as the right hand edge $(\XGTd^1_1(t), \XGTd^2_2(t), \ldots, \XGTd^n_n(t); t \in \Nnint)$ when $\XGTd$ has initial distribution $M_z$, $z = (0,\ldots,0)$ (i.e. 
$\XGTd^j_k(0) = 0$, $1 \leq k \leq j \leq n$).
But, a version of the right hand edge may be constructed from paths of 
$Z$ begun at the origin so that Theorem \ref{t:discretegeometricintro} gives

\begin{equation}\label{e:johanssonident}
 \left( Z_n^\dagger(t); t \geq 0\right) \overset{\mathrm{dist}}{=} \left( \max_{\pi \in \Pi(t,n)} \sum_{(i,j) \in \pi} (Z_j(i) - Z_j(i-1)); t \geq 0\right).
\end{equation}

As a consequence, Theorem \ref{t:discretegeometricintro}
provides a new proof that such last passage percolation times have the same distribution as 
the rightmost particle in the conditioned process $Z^\dagger$ (the distribution of which, at a fixed time, is given by the Meixner ensemble -- see Johansson \cite{MR1737991} or \cite{johansson2007mdm}).  This is a key step in obtaining the Tracy-Widom distribution
in this setting.

Note that the dynamics discussed above are different to those exhibited in \cite{borodin2008agr} for geometric jumps.  In particular, the particles in
the process we described above are blocked by the position of the particle
immediately above and to the right of them at the \emph{previous} time step.

\subsection{With wall at the origin}
The final example of the paper uses the ideas introduced above to 
construct a continuous time process $(\XGTw(t);\; t \geq 0)$ on a \emph{symplectic} Gelfand-Tsetlin cone.  
The latter are so termed because they are in direct correspondence with the symplectic tableau 
arising from the representations of the symplectic group \cite{sundaram90}.

The space $\SympKone_{n}$ of integer valued symplectic Gelfand-Tsetlin patterns 
may be defined (see for example \cite{MR2427082} or \cite{proctor1988osg}) as 
the set of point configurations 
$(x^1, x^2, \ldots, x^{n-1}, x^{n})$ such that 
\begin{itemize}
 \item $x^{2i-1}, x^{2i} \in \Nnint^{i}$ for $1 \leq i \leq \lfloor \frac{n}{2} \rfloor$ and 
$x^{n} \in \Nnint^{(n+1)/2}$ if $n$ is odd,
 \item $x^{2i-1} \prec x^{2i}$ for $1 \leq i \leq \lfloor \frac{n}{2} \rfloor$,
 \item $x^{2i} \preceq x^{2i+1}$ for $1 \leq i \leq \lfloor \frac{n-1}{2} \rfloor$.
\end{itemize}

So the all the points in a symplectic pattern lie to the right of an impenetrable wall at the origin, 
represented diagrammatically below.

\begin{center}
\begin{tabular}{|*{9}{c}}
 $x_1^1$ &	& 	& 	&&	&	&& \\
         &$x^2_1$&	&	&&	&& \\
 $x^3_1$&	&$x_3^2$&&	&&& \\
 	 &$x^4_1$&	&$x^4_2$&	&&& \\
 $x^5_1$&	&$x^5_2$&	&$x^5_3$&& \\

 &$\vdots$&&	$\vdots$& &$\ddots$& \\
 $x^n_1$& & $x^n_2$& & $x^n_3$ &$\ldots$& $x^n_n$ \\
\end{tabular}
\end{center}

In the vein of previous sections, we construct a process 
$(\XGTw(t),t \geq 0)$ in 
$\SympKone_{n}$ 
in which only one particle jumps of its own volition at each instant and 
a blocking and pushing interaction maintains the interlacing constraints.

Fix $q \in (0,1)^n$. The top particle $\XGTw^1_1$ jumps right at rate $q_1$ 
and left at rate $q^{-1}_1$, \emph{apart} from at origin where its left jumps are suppressed.  The second row also only has one particle, $\XGTw^1_2$, which 
jumps rightwards at rate $q_1^{-1}$ and leftwards at rate $q_1$ (notice rates are reversed), \emph{except} 
at instances when $\XGTw^1_1(t) = \XGTw^2_1(t)$.  In the latter case, it is pushed 
rightwards if $\XGTw^1_1$ jumps to the right and any leftward jumps are suppressed.

The remaining particles evolve in a similar fashion -- 
on row $2k - 1$, particles take steps to the right at rate $q_k$ and left 
at rate $q^{-1}_k$ when they are not subject to the blocking or pushing
required to keep the process in the state space, in particular 
$\XGTw^{2k-1}_1$ has any leftward jump from the origin suppressed.
On row $2k$, the rates are reversed but the same blocking and pushing mantra applies.

We will deduce that for appropriate initial conditions, the marginal 
distribution of each row $(\XGTw^k(t); t \geq 0)$ is a 
Markov process.  The $Q$-matrices for the marginal processes can be written in terms of
\emph{symplectic} Schur functions, the definition of which is 
similar to that of the classic Schur function \eqref{e:schur} --
 they are sums over geometrically weighted symplectic Gelfand-Tsetlin
patterns.

Fix $k > 0$ and suppose that either $n = 2k-1$ or $n = 2k$. 
Now let $\Sxwk^k = \{ z \in \Z^k : 0 \leq z_1 \leq z_2 \ldots \leq z_k \}$ 
and define $\SympKone_{n}(z)$ to be the set of symplectic patterns 
$x$ in $\SympKone_{n}$ with bottom row $x^{n}$ equal to $z \in \Sxwk^k$.  
The geometric weight $w_n^q$ on $\SympKone_n$ is
\[
 w_{2k-1}^{(q_1,\ldots, q_k)}(x) = q_k^{|x^{2k-1}| - |x^{2k-2}|} \prod_{i=1}^{k-1} q_i^{|x^{2i-1}| - |x^{2i-2}|+|x^{2i-1}|-|x^{2i}|}
\]
and
\[
 w_{2k}^{(q_1, \ldots, q_k)}(x)= \prod_{i=1}^k q_i^{|x^{2i-1}| - |x^{2i-2}|+|x^{2i-1}|-|x^{2i}|},
\]
using the convention that $|x^0| = 0$ and empty products are equal to 1 (so 
$w_1^{(q_1)}(x) = q_1^{|x^1|}$).

Then, the symplectic Schur function $Sp^{n}_z:\R^k \to \R$, $z \in \Sxwk^k$, 
$k \geq 1$ is defined (see \cite{MR1426737}) by
\begin{equation}\label{e:oddsympschur}
 Sp^{n}_z(q_1, \ldots, q_k) = \sum_{x \in \SympKone_{n}(z)} w_{n}^q(x).
\end{equation}

For even $n$, $Sp^{n}$ gives the characters of irreducible representations of 
the symplectic group $Sp(n)$ \cite{sundaram90}.  For odd $n$, $Sp^n$ 
was introduced by Proctor \cite{proctor1988osg} and can interpretted as the character
of the irreducible representations of a group that interpolates between the classical 
groups $Sp(n)$ and $Sp(n+1)$ \cite{Maliakas1999640}.

Define the $Q$-matrix $Q_{n}:\Sxwk^k \times \Sxwk^k \to \R$ as follows.
For $x \in \Sxwk^k$ and $x \pm e_i \in \Sxwk^k$, some $1\leq i \leq k$,
\begin{equation}\label{e:Q2km1}
 Q_{n}(x,x \pm e_i) = \frac{Sp^{n}_{x \pm e_i}(q)}{Sp^{n}_x(q)}.
\end{equation}

All other off diagonal entries vanish and the diagonals are given by
\begin{equation}
 -Q_{2k-1}(x,x) = \sum_{i=1}^{k-1} \left( q_i + q_i^{-1}\right) + q_k^{-1} \Indi_{[x_1 > 0]} + q_k,
\end{equation}
and
\begin{equation}\label{e:Q2k}
 -Q_{2k}(x,x) = \sum_{i=1}^k \left( q_i + q_i^{-1}\right).
\end{equation}

A corollary of the intertwinings we prove in sections \ref{s:poissonwallI} and \ref{s:poissonwallII} is that $Q_n$ is conservative.

For $z \in \Sxwk^k$ define $M^{n}_z:\SympKone_n \to [0,1]$ by
\[
 M^{n}_z(x) = \frac{w^q_{n}(x)}{Sp^{n}_z(q)}.
\]
Then the definitions of the symplectic Schur functions imply that $M^n_z(\cdot)$ gives a probability distribution on patterns in $\SympKone_n(z)$.  

From these ingredients we obtain

\begin{theorem}\label{t:wallevenX}
Suppose $\XGTw$ has initial distribution given by $M^n_z(\cdot)$, then
$(\XGTw^{n}(t); t \geq 0)$ is distributed as a Markov process 
with $Q$-matrix $Q_n$, started from $z$.
\end{theorem}

The relevance of this theorem to the discussion in the introduction may 
again be seen by examining the evolution of the right hand edge of $\XGTw$.
Suppose we have a system of $n$ particles with positions
$(\XGTw^1_1(t),\XGTw^2_1(t)+1, \XGTw^3_2(t)+2, \ldots, \XGTw_{\lfloor (n+1)/2 \rfloor}^n(t)+n-1; t \geq 0)$.  

Particle $i > 1$ attempts to jump rightwards at rate $\gamma_i = q_{(i+1)/2}$ if $i$ is odd or $\gamma_i = q^{-1}_{i/2}$ if $i$ is even and leftwards at rate $\gamma^{-1}_i$.  An attempted left jump succeeds only if the destination site is 
vacant, otherwise it is suppressed.  A rightward jump always succeeds, and, any
particle occupying the destination site is pushed rightwards.  A particle being pushed rightwards also pushes any particle standing in its way, so a rightward jump by a particle could cause many particles to be pushed.
So far we have essentially described the dynamics of the ``PushASEP'' process introduced in \cite{MR2438811}.  Our process differs by the presence of a wall: the leftmost particle (identified with $\XGTw^1_1$) is modified so that any leftward jump at the origin suppressed.  Also, the particle rates are restricted in that for odd $i$, the jump rates of particle $i$ and $i+1$ are inverses of each other (which is not the case in \cite{MR2438811}).

As in the previous examples, the bottom row $(\XGTw^n(t); t \geq 0)$ 
may be realised as a Doob $h$-transform and we deduce identities analogous to \eqref{e:oconnellyorident} and \eqref{e:johanssonident}.  For simplicity, 
we shall only consider the case that $n = 2k$.  The case of odd $n$ can be treated 
with similar arguments but it is complicated slightly due to the non-standard behaviour of 
$\XGTw_1^n$ at the wall.

Let $Z$ be a $k$-dimensional random walk
in which the $i^{th}$ component jumps rightwards at rate $q^{-1}_i$ and leftwards
at rate $q_i$.  It is readily seen that $Q_{2k}$ is the $Q$-matrix of $Z^\dagger$, the $h$-transform of $Z$ killed on leaving $\Sxwk^k$ under harmonic functions
\[
h_{2k}(x) = q_1^{-x_1} q_2^{-x_2} \ldots q_k^{-x_k} Sp^{2k}_x(q),\; x \in \Sxwk^k.
\]

Theorem \ref{t:wallevenX} shows that 
$\left( Z^\dagger_k(t); t \geq 0 \right)$ has the same law as $\left(\XGTw^{2k}_k(t); t \geq 0\right)$ when $\XGTw$ is initially distributed according to $M^{2k}_z$ and $Z(0) = z \in \Sxwk^k$.

But if $z = (0,\ldots,0)$, a process with
the same law as the right hand edge of $\XGTw$ can be constructed from the paths of $Z$ and a random walk $\tilde Z$ that is independent of, but identically
distributed to $-Z$.  The resulting identity in distribution can be stated succinctly 
in terms of the $2k$-dimensional random walk $\bar Z(t) = (Z_1(t),\tilde Z_1(t), Z_2(t),\tilde Z_2(t), \ldots, Z_k(t),\tilde Z_k(t))$ as follows

\[
\left( Z^\dagger_k(t); t \geq 0 \right) 
 \overset{dist}{=} \left( \sup_{0 \leq t_1 \leq \ldots \leq t_{2k+1}=t} \sum_{i=1}^{2k} \left(\bar Z_i(t_{i+1}) - \bar Z_i(t_{i}) \right); t \geq 0 \right).
\]

The Brownian analogue of this result will be considered in 
\cite{BFPSWmaxdysonbmboundary09}.

\section{Proof of Theorem \ref{t:charlierresult}}\label{s:poisson}

Let $( \XGT(t);\; t \geq 0)$ be the process on $\Kone_{n}$ satisfying
 the dynamics described in section \ref{s:poissonintro}.  
It is clear from this 
description that the law of
$\{ \XGT^{n}(t);\; t \geq 0\}$ is conditionally independent
of $\{ \XGT^{j}(t);\; t \geq 0, j < n - 1\}$ 
given $\{ \XGT^{n - 1}(t);\; t \geq 0\}$.  That is, the dynamics of the particle in row $n$ depend on the evolution of particles in the rows above only through the particles in row $n-1$.  Hence the theorem may be proven inductively by studying only the bottom and
penultimate layers of the pattern.

To this end, we assume for induction that the conclusion of 
\ref{t:charlierresult} holds. Then, when $\XGT(0)$ is distributed according to 
$M_z(\cdot)$, the bottom layer $(\XGT^{n}(t); t \geq 0)$ is Markovian and 
evolves according to the conservative $Q$-matrix  $Q_X$ defined via
\[
Q_X(x,x+e_i) = \frac{S_{x + e_i}(q)}{S_{x}(q)}\Indi_{[x+e_i \in \Sx]}, \quad 1 \leq i \leq n, \; x \in \Sx,
\]
and all other off diagonal entries set to zero.

We will define a Markov process $(X(t), Y(t); t\geq 0)$ on $\Ss = \{ (x,y) \in \Sx \times \Sy : x \preceq y\}$ (recall $x \preceq y$ means that $y_i \leq x_i \leq y_{i+1}$, $1 \leq i \leq n$)  in which $X$ evolves according to $Q_X$ while $Y$ 
evolves independently of $X$ apart from the blocking and pushing interaction.
One should think of $(X,Y)$ as the penultimate and bottom layer of our construction
in $\Kone_{n+1}$.  So, to complete the induction step it is sufficient to show that marginally $Y$ is Markovian and evolves according to
\[
  Q_Y(y,y+e_i) = q_i \frac{\tilde h(y + e_i)}{\tilde h(y)},
\]
for $y \in \Sy$, where $\tilde h$ is given by
\[
  \tilde h(y) = q_1^{-x_1} \ldots q_{n+1}^{-x_{n+1}} S_y(q_1, \ldots, q_{n+1}),\; y \in \Sy.
\]
for some $q_{n+1} > 0$, and all other off diagonal entries vanish.  The diagonal entries are given by
\[
 Q_Y(y,y) = -\sum_{i=1}^{n+1} q_i.
\]

Appropriate dynamics for $(X,Y)$ are specified by the conservative $Q$-matrix $\Gen$ with off diagonal entries given by
\[
\Gen((x,y),(x^\prime,y^\prime)) = \left\{\begin{array}{ll}
Q_X(x,x + e_i), &\; (x^\prime,y^\prime) = (x + e_i,y), \; x_i < y_{i+1},  \\
Q_X(x,x + e_i), &\; (x^\prime,y^\prime) = (x + e_i,y + e_{i+1}), \; x_i = y_{i+1}, \\
q_{n+1}, &\; (x^\prime,y^\prime) = (x,y + e_j), \\
0 &\; \mathrm{otherwise}
\end{array}\right.
\]
for $(x,y),(x^\prime,y^\prime) \in \Ss$ and $1 \leq i \leq n$, $1 \leq j \leq n+1$.  
The diagonal entry $-\Gen((x^\prime,y^\prime),(x^\prime,y^\prime))$ is given by
\begin{equation}\label{e:Arateofleaving}
\sum_{i=1}^n q_i + q_{n+1} \sum_{i=1}^n \Indi_{[y^\prime_i < x^\prime_i]} +  q_{n+1} 
= \sum_{i=1}^{n+1} q_i + q_{n+1} \sum_{i=1}^n \Indi_{[y^\prime_i < x^\prime_i]}. 
\end{equation}

Now, as an immediate consequence of the definition of the Schur function in \eqref{e:schur}, we have
\begin{eqnarray*}
 S_{z}(q_1,q_2, \ldots, q_n, q_{n+1}) &=& \sum_{x \in \Kone_n(z)} w^q(x) \\
& = & \sum_{z^\prime \in \Sx : z^\prime \preceq z} q_{n+1}^{|z|-|z^\prime|} \sum_{x \in \Kone_n(z^\prime)} \prod_{j=1}^{n} q_j^{|x^j| - |x^{j-1}|} \\
& = & \sum_{z^\prime \preceq z} q_{n+1}^{|z| - |z^\prime|} S_{z^\prime}(q_1, \ldots, q_n).
\end{eqnarray*}

So the marginal distribution of the penultimate row of particles under 
the initial distribution defined in \eqref{e:distonKn} is given by 
$m(\cdot,y)$ where $y \in \Sy$ is fixed and $m:\Ss \to [0,1]$ is defined by
\[
m(x,y) = q_{n+1}^{|y| - |x|}\frac{S_x(q)}{S_y(\tilde q)},
\]
where $\tilde q = (q_1,q_2, \ldots, q_n, q_{n+1})$.

Furthermore,
\begin{equation}\label{e:kernellambda}
 \Lambda(y,(x^\prime, y^\prime)) = m(x^\prime, y^\prime) \Indi_{[y^\prime = y]}.
\end{equation}
defines a Markov kernel from $\Sy$ to $\Ss$.  That is, for each $y \in \Sy$, 
$\Lambda(y,\cdot)$ defines a probability distribution on $\Ss$.

The heart of our proof is showing that the conservative $Q_Y$ is \emph{intertwined} with $\Gen$ via $\Lambda$,
\begin{equation}\label{e:intqmatrix}
 Q_Y \Lambda = \Lambda \Gen.
\end{equation}

From here, lemma \ref{l:qmatintimpkernint} shows that $\Lambda$ intertwines the corresponding transition kernels.  That is, if $(p_t; t \geq 0)$ are the transition kernels corresponding to $Q_Y$ and $(q_t; t \geq 0)$ those to $\Gen$,
then for $y \in  \Sy$, $(x^\prime, y^\prime) \in \Ss$ and $t \geq 0$,
\[
 p_t(y, y^\prime)m(x^\prime,y^\prime) = \sum_{x \prec y} m(x,y) q_t((x,y), (x^\prime, y^\prime)),
\]

An immediate consequence of this relationship is that for bounded $f:\Sy\to\R$,
\begin{eqnarray*}
 \E^y[f(Y(t))] & =& \sum_{(x^\prime, y^\prime)} \sum_{x \prec y} m(x,y) q_t((x,y), (x^\prime, y^\prime)) f(y^\prime) \\
& = & \sum_{(x^\prime, y^\prime)} p_t(y, y^\prime) m(x^\prime,y^\prime) f(y^\prime)\\
& = & \sum_{y^\prime} p_t(y,y^\prime) f(y^\prime) \sum_{x^\prime \prec y^\prime} m(x^\prime, y^\prime)\\
& = & \sum_{y^\prime} p_t(y,y^\prime) f(y^\prime),
\end{eqnarray*}
where $\E^y$ is the expectation operator corresponding to the measure under which $(X,Y)$ 
has initial distribution $\Lambda(y,\cdot)$.

When $0 \leq t_1 \leq \ldots \leq t_N$ and $f_1, \ldots, f_N:\Sy\to\R$ are bounded, 
the preceeding argument generalises and the intertwining shows that
\[
 \E^y[f_1(Y(t_1)) \ldots f_N(Y(t_N))] = \sum_{y^1, \ldots, y^N} p_t(y,y^1) p_t(y^1,y^2) \ldots p_t(y^{N-1}, y^N) f_1(y^1) \ldots f_N(y^N).'
\]

This is essentially the argument of Rogers and Pitman \cite{rogers1981mf} and establishes

\begin{theorem}\label{t:poisson}
Suppose $(X(t),Y(t); t \geq 0)$ is a Markov process with $Q$-matrix 
$\Gen$ and initial distribution $\Lambda(y,\cdot)$, for some $y \in \Sy$.
Then $Q_Y$ and $\Gen$ are interwined via $\Lambda$ and as a consequence, 
$(Y(t); t \geq 0)$ is distributed as a Markov process with $Q$-matrix $Q_Y$, 
started from $y$.
\end{theorem}

The intertwining \eqref{e:intqmatrix} is equivalent to
\begin{equation}\label{e:QYintertwining}
Q_Y(y,y^\prime) = \sum_{x \preceq y} \frac{m(x,y)}{m(x^\prime,y^\prime)} \Gen((x, y),(x^\prime, y^\prime)), \;\; y \in \Sy, \; (x^\prime, y^\prime) \in \Ss,
\end{equation}
where the summation is over the points $x$ in $\Sx$ that interlace with $y$. 
As the particles can only make unit jumps rightwards, both sides of the 
expression vanish unless either $y^\prime = y$ or $y^\prime = y + e_j$, 
for some $1 \leq j \leq n+1$.

We first consider the case when $y = y^\prime$, corresponding to the diagonal entries of
$Q_Y$.  The right hand side of the expression is
\[
 \sum_{x \preceq y^\prime} \frac{m(x,y^\prime)}{m(x^\prime,y^\prime)} \Gen((x, y^\prime),(x^\prime, y^\prime)).
\]
Using the definition of $m$, this becomes
\begin{equation}\label{e:QYintertwiningdiag}
 \sum_{x \preceq y^\prime} q_{n+1}^{|x^\prime| - |x|} \frac{S_{x}(q)}{S_{x^\prime}(q)} \Gen((x, y^\prime),(x^\prime, y^\prime)).
\end{equation}
Now, $\Gen((x, y^\prime),(x^\prime, y^\prime))$ is non zero for $x \preceq y^\prime$ only if 
$x = x^\prime$ or $x = x^\prime - e_i$ for some $1 \leq i \leq n$.  When $x = x^\prime$, $-\Gen((x, y^\prime),(x^\prime, y^\prime))$ is the rate of
leaving at $(x^\prime, y^\prime)$, given in \eqref{e:Arateofleaving}.  On the other hand if $x = x^\prime - e_i$, $\Gen((x, y^\prime),(x^\prime, y^\prime))$ 
is the rate at which the $i^{th}$ $X$ particle jumps rightwards (without pushing
a $Y$ particle).  But, such values of $x$ are included in the summation only if 
$x = x^\prime - e_i \preceq y = y^\prime$, i.e. $x_i^\prime > y_i^\prime$.

Combining this with \eqref{e:QYintertwiningdiag} and \eqref{e:Arateofleaving}
and the fact that $q_{n+1}^{|x^\prime| - |x^\prime - e_i|} = q_{n+1}$,
we see that if $y = y^\prime$ the right hand side of \eqref{e:QYintertwining} is
\[
\sum_{i=1}^n q_{n+1} \frac{S_{x^\prime - e_i}(q)}{S_{x^\prime}(q)} Q_X(x^\prime - e_i,x^\prime)\Indi_{[x_i^\prime > y_i^\prime ]} - \sum_{i=1}^{n+1} q_i -  \sum_{i=1}^{n} q_{n+1}\Indi_{[y_i^\prime < x_i^\prime]}.
\]

The first summand above is
\[
 q_{n+1} \frac{S_{x^\prime - e_i}(q)}{S_{x^\prime}(q)} Q_X(x^\prime - e_i,x^\prime)\Indi_{[x_i^\prime > y_i^\prime]} = q_{n+1} \Indi_{[y_i^\prime < x_i^\prime]},
\]
so the first and last summations above disappear and we are left with $-\sum_{i=1}^{n+1} q_i$, which is exactly $Q_Y(y^\prime, y^\prime)$.

If $y \neq y^\prime$, the only other possibility is that $y^\prime = y + e_i$ for some 
$1 \leq i \leq n+1$.  Let us first deal with the simplest case, where $i = 1$, that is, $y^\prime = y + e_1$. 
The only value of $x$ for which $\Gen( (x , y^\prime - e_1), (x^\prime, y^\prime) )$ is non zero is $x = x^\prime$ as the first $Y$ particle is never pushed by an $X$ particle.  Furthermore, $y^\prime_1 - 1 < y^\prime_1 \leq x^\prime_1$ and so the jump of $Y_1$ is certainly not blocked.
Hence, 
\begin{eqnarray*}
\sum_{x : x \preceq y} \frac{m(x,y)}{m(x^\prime,y^\prime)} \Gen((x, y),(x^\prime, y^\prime))
 &=& \frac{m(x^\prime,y^\prime - e_1)}{m(x^\prime,y^\prime)} \Gen((x^\prime,y^\prime - e_1),(x^\prime, y^\prime))\\
& = & \frac{q_{n+1}^{|y^\prime - e_1| - |x^\prime|}}{q_{n+1}^{|y^\prime| - |x^\prime|}}
\frac{S_{x^\prime}(q)}{S_{x^\prime}(q)}
\frac{S_{y^\prime}(\tilde q)}{S_{y^\prime - e_1}(\tilde q)} q_{n+1}\\
& = & Q_Y(y^\prime - e_1,y^\prime).
\end{eqnarray*}

So in this case, \eqref{e:QYintertwining} is satisfied.

For $i > 1$, consider the dichotomy $x^\prime_{i-1} < y^\prime_i$
or $x^\prime_{i-1} = y^\prime_i$.  Suppose we are in the former case, i.e. $y^\prime = y + e_i$ and $x^\prime_{i-1} < y^\prime_i$. 
It is not possible that the movement in the $i^{th}$ component of $Y$ could have 
been instigated due to pushing by the $(i-1)^{th}$ $X$ particle (a push could only have occurred if $x^\prime_{i-1} - 1 = y^\prime_i - 1$).  
Thus, as in the $i = 1$ case above, $\Gen( (x , y^\prime - e_i), (x^\prime, y^\prime) )$ is non zero only for $x = x^\prime$ and almost identical calculations verify \eqref{e:QYintertwining}.

The second $i > 1$ subcase is that $x^\prime_{i-1} = y^\prime_i$ and $y = y^\prime - e_i$.  
Here the only possibility is that the $i^{th}$ $Y$ particle 
``did not jump but was pushed'',
which one may confirm by noting that $x^\prime$ does not interlace with $y^\prime - e_i$
when $x^\prime_{i-1} = y^\prime_i$.  So, the right hand side of \eqref{e:QYintertwining} is given by
\[
 \frac{m(x^\prime-e_{i-1}, y^\prime-e_{i})}{m(x^\prime,y^\prime)}\Gen((x^\prime-e_{i-1},y^\prime-e_{i}),(x^\prime,y^\prime)).
\]
Using the definitions of $m$ and $\Gen$, this becomes
\[
 \frac{q_{n+1}^{|y^\prime-e_{i}| - |x^\prime-e_{i-1}|}}{q_{n+1}^{|y^\prime| - |x^\prime|}} \frac{S_{x^\prime-e_{i-1}}(q)}{S_{y^\prime-e_{i}}(\tilde q)}\frac{S_{y^\prime}(\tilde q)}{S_{x^\prime}(q)}\frac{S_{x^\prime}(q)}{S_{x^\prime-e_{i-1}}(q)},
\]
a quantity which is easily seen to equal $Q_Y(y^\prime - e_i, y^\prime)$.

This concludes the proof that $Q_Y$ and $\Gen$ are intertwined via $\Lambda$.

\section{Proof of Theorem \ref{t:discretegeometricintro}}\label{s:discretegeometric}

It is again sufficient to consider any pair of consecutive rows $(X,Y)$ and 
construct the process iteratively.

Let $(X(t); t \in \Nnint)$ be an $n$ dimensional Markov chain in $\Sx$ with one step transition kernel
\[
 p_X(x, x^\prime) = a(q) \frac{S_{x^\prime}(q)}{S_x(q)} \Indi_{[x \prec x^\prime]}.
\]
where $q \in (0,1)^n$, $a(q)= \prod_{i=1}^n (1-q_i)$ and for $x,x^\prime\in \R^n$,
$x \prec x^\prime$ indicates that the inequality $x_1 \leq x_1^\prime \leq x_2 \leq \ldots \leq x^\prime_{n-1} \leq x_n \leq x_n^\prime$ holds.

Let $\xi_i(t)$ ($t\in \Nnint, 1 \leq i \leq n+1$) be geometric($q_{n+1}$) random variables that are independent of each other and of $X$,
\[
 \Prob(\xi_i(k) = j) = (1 - q_{n+1})q_{n+1}^j, \quad j = 0,1,2,\ldots.
\]

Define a process $(Y(t); t \in \Nnint)$ in $\Sy$ in terms of $X$ using the recursion
\[
  Y_1(t+1) = \min(Y_1(t) + \xi_1(t+1), X_1(t)),
\]
\[
 Y_{n+1}(t+1) =  \max(Y_{n+1}(t),X_n(t+1)) + \xi_{n+1}(t+1), 
\]
\[
 Y_j(t+1) = \min(\max(Y_j(t), X_{j-1}(t+1)) + \xi_j(t+1), X_j(t)),
\]
for $2 \leq j \leq n$.

The recursion encodes 
the blocking and pushing mechanism, maintaining the initial interlacing relationship, so $X(t) \prec Y(t)$ for each $t$.

We will prove that if $\Lambda$ is as defined in \eqref{e:kernellambda} then

\begin{theorem}\label{t:discretegeometric}
If $(X,Y)$ is initially distributed according to $\Lambda(y,\cdot)$, $y \in \Sy$,
and then evolves according to the recursion above, 
the marginal process $(Y(t); t \geq 0)$ is distributed as an 
$n+1$ dimensional Markov process with transition kernel
\[
 p_Y(y, y^\prime) = a(q_1,\ldots,q_{n+1}) \frac{S_{y^\prime}(q_1,\ldots,q_{n+1})}{S_y(q_1,\ldots,q_{n+1})} \Indi_{[y \prec y^\prime]}.
\]
started from $y$.
\end{theorem}

Our strategy, again, is to prove that $\Lambda$ interwines the corresponding transition probabilities.  Suppose $(x,y), (x^\prime, y^\prime) \in \Ss$, $x \prec x^\prime$ and $y \prec y^\prime$.  
Let us write down $q((x,y), (x^\prime,y^\prime))$, the one step transition probabilities for $(X,Y)$.  Firstly note that
\[
q((x,y), (x^\prime,y^\prime)) = 
r(y^\prime,x^\prime,x,y)p_X(x,x^\prime).
\]
where
\[
r(y^\prime,x^\prime,x,y) = \Prob(Y(1)=y^\prime|X(1)=x^\prime, X(0)=x, Y(0)=y).
\]

Using the definition of $Y$, $r$ can be conveniently expressed in terms of 
the ``blocking'' and ``pushing'' factors $b,c:\Z^2 \to \R_+$
\[
 b(u,v) = (1-q_{n+1})\Indi_{[v < u]} + \Indi_{[u = v]},
\]
\[
c(u,v) = q_{n+1}^{-v}\Indi_{[u \leq v]} + q_{n+1}^{-u}\Indi_{[u > v]}.
\]
Then $r(y^\prime,x^\prime,x,y)$ is equal to
\[
 q_{n+1}^{y_1^\prime - y_1} b(x_1, y_1^\prime) 
\left(\prod_{i=2}^{n} q_{n+1}^{y_i^\prime} b(x_{i}, y_{i}^\prime) c(x_{i-1}^\prime,y_i)\right) q_{n+1}^{y_{n+1}^\prime}(1-q_{n+1})c(x_n^\prime, y_{n+1}).
\]

To prove the theorem we will need the following ``integrating out'' lemma.  

\begin{lemma}\label{l:bpsum}
Suppose $v_2, v_1^\prime, u^\prime \in \Z$ satisfy $v_1^\prime \leq v_2$, $v_1^\prime \leq u^\prime$ and $q_{n+1} \neq 1$.  Then we have 
\begin{equation}
\sum_{u = v_1^\prime}^{v_2 \wedge u^\prime} q_{n+1}^{-u} b(u, v_1^\prime)c(u^\prime,v_2) = q_{n+1}^{-u^\prime - v_2}
\end{equation}
\end{lemma}

The lemma may be understood more readily by imagining that we 
are considering the $n=1$ case,
so that there is one ``$X$'' particle nestled between two ``$Y$'' particles.  
We may fix the initial and final positions of the ``$Y$'' particles ($v$ and $v^\prime$ in the lemma above) 
and also the final position of the ``$X$'' particle ($u$ in the lemma) -- it is the starting location 
of the $X$ particle that we are integrating out.  The summation is over the possible values that the $X$ particle \emph{may} have started from.  
It must be at least equal to the final position of the left most $Y$ particle $v_1^\prime$, as this particle cannot 
overtake the $X$ particle (see recursion equations above). Also, it cannot exceed either the initial position 
of the second $Y$ particle $v_2$ (due to the interlacing constraint) or the final position of the $X$ particle $u^\prime$ (as the particles may only jump rightwards).

\begin{proof}

After using the definitions of $b$ and $c$, the sum becomes
\[
 \sum_{u = v_1^\prime}^{v_2 \wedge u^\prime} q_{n+1}^{-u}  \left( (1-q_{n+1})\Indi_{[v_1^\prime < u]} + \Indi_{[u = v_1^\prime]} \right) \left( q_{n+1}^{-v_2}\Indi_{[u^\prime \leq v_2]} + q_{n+1}^{-u^\prime}\Indi_{[u^\prime > v_2]}\right)
\]

Now expand the brackets in the summand and sum the terms individually.
We find
\[
 \sum_{u = v_1^\prime}^{v_2 \wedge u^\prime}(1-q_{n+1})\Indi_{[v_1^\prime < u]} q_{n+1}^{-v_2-u}\Indi_{[u^\prime \leq v_2]} = q_{n+1}^{-v_2}\left(q_{n+1}^{-u^\prime} - q_{n+1}^{-v_1^\prime} \right)\Indi_{[u^\prime \leq v_2]},
\]

\[
\sum_{u = v_1^\prime}^{v_2 \wedge u^\prime}(1-q_{n+1})\Indi_{[v_1^\prime < u]}q_{n+1}^{-u^\prime-u}\Indi_{[u^\prime > v_2]} = q_{n+1}^{-u^\prime}\left(q^{-v_2}_{n+1} - q^{-v^\prime_1}_{n+1}\right)\Indi_{[u^\prime > v_2]},
\]
\[
 \sum_{u = v_1^\prime}^{v_2 \wedge u^\prime} \Indi_{[u = v_1^\prime]} q_{n+1}^{-v_2 - u}\Indi_{[u^\prime \leq v_2]}  = q_{n+1}^{-v_2 - v_1^\prime}\Indi_{[u^\prime \leq v_2]}
\]
and
\[
\sum_{u = v_1^\prime}^{v_2 \wedge u^\prime} \Indi_{[u = v_1^\prime]} q_{n+1}^{-u^\prime - u}\Indi_{[u^\prime > v_2]}= q_{n+1}^{-u^\prime - v_1^\prime}\Indi_{[u^\prime > v_2]}.
\]
Summing the above expressions gives the result.
\end{proof}

The interesting thing about this scheme, as we will see in a moment, is that we may apply it successively
from left to right when there are $n$ particles so that the leftmost 
particles get heavier and heavier until we have reduced the problem 
to the $n=1$ case.

When the initial distribution is $\Lambda(y,\cdot)$, the joint distribution after one time step is given by
\[
\pi(x^\prime, y^\prime) = \sum_{x \prec y} m(x,y) q((x,y), (x^\prime,y^\prime)).
\]

Expanding the sum and incorporating the conditions $y_i^\prime \leq x_i$ 
and $x \prec x^\prime$ into the summation indices yields
\begin{equation}\label{e:rhs}
\pi(x^\prime, y^\prime) = \sum_{x_n = y_n^\prime}^{y_{n+1} \wedge x_n^\prime}\ldots\sum_{x_1 = y_1^\prime}^{y_2\wedge x_1^\prime} m(x,y) q((x,y), (x^\prime,y^\prime)).
\end{equation}

The summand in \eqref{e:rhs} equals
\begin{eqnarray*}
& & q_{n+1}^{y_1^\prime + y_2 + \ldots + y_{n+1} + y_{n+1}^\prime - x_1 - \ldots -x_n}
 b(x_1, y_1^\prime)\left(\prod_{i=2}^{n} q_{n+1}^{y_i^\prime} b(x_{i}, y_{i}^\prime)c(x_{i-1}^\prime,y_i)\right) \\
& & \quad \times c(x_n^\prime, y_{n+1})a(\tilde q) \frac{S_{x^\prime}(q)}{S_y(\tilde q)},
\end{eqnarray*}
for $x \preceq y, x^\prime \preceq y^\prime, x \prec x^\prime, y \prec y^\prime$ and vanishes elsewhere.

Now, one notices that we may use lemma \ref{l:bpsum} to iteratively evaluate the summation over $x_1, x_2, \ldots, x_{n}$ (in that order).  
More concretely, first apply the lemma with $u^\prime = x_1^\prime, v = (y_1, y_2), v^\prime = (y_1^\prime, y_2^\prime)$ to reveal that the sum $\sum_{x_1 = y_1^\prime}^{y_2 \wedge x_1^\prime} m(x,y) Q((x,y), (x^\prime,y^\prime))$ is equal to

\begin{eqnarray*}
&&	q_{n+1}^{y_1^\prime + y_2^\prime + y_3 + \ldots + y_{n+1} + y_{n+1}^\prime - x_1^\prime -x_2 - \ldots -x_n} b(x_2, y_2^\prime) \\
&& \quad \times \left(\prod_{i=3}^{n} q_{n+1}^{y_i^\prime} b(x_{i}, y_{i}^\prime) c(x_{i-1}^\prime,y_i)\right) c(x_n^\prime, y_{n+1})a(\tilde q) \frac{S_{x^\prime}(q)}{S_y(\tilde q)}.
\end{eqnarray*}

This expression is again in a suitable form to apply lemma \ref{l:bpsum}, but this time with $u^\prime = x_2^\prime, v = (y_2, y_3), v^\prime = (y_2^\prime, y_3^\prime)$ and summing over $x_2$.  Continuing in this 
fashion shows that \eqref{e:rhs} is equal to 
\begin{equation}
q_{n+1}^{y_1^\prime + y_2^\prime \ldots + y_{n+1}^\prime - x_1^\prime  - \ldots - x_{n}^\prime}a(\tilde q) \frac{S_{x^\prime}(q)}{S_y(\tilde q)} = m(x^\prime,y^\prime) p_Y(y, y^\prime).
\end{equation}

Hence we have verified the intertwining
\[
m(x^\prime,y^\prime) p_Y(y, y^\prime) = \sum_{x \preceq y} m(x,y) q((x,y), (x^\prime,y^\prime)),
\]
and Theorem \ref{t:discretegeometric} follows from the argument of \cite{rogers1981mf} discussed in the previous section.

\section{Proof of Theorem \ref{t:wallevenX}}

As in the previous two examples, we give a row by row construction.  This time
the asymmetry between odd rows and even rows means we have to specify how to iterate from even rows to odd rows
and odd rows to even rows separately (presented below in \ref{s:poissonwallI} 
and \ref{s:poissonwallII} respectively).

En route to proving Theorem \ref{t:wallevenX}, we need to conclude that 
$Q_{n}$ is a conservative $Q$-matrix for each $n$.

This will be achieved by an inductive argument.  Let H($n$) denote the hypothesis that 
$Q_n$ is a conservative $Q$-matrix. It is easy to establish H(1),
that $Q_1$ is conservative -- recall that for $x_1 \geq 0$, $Sp^1_{(x_1)} = q_1^{x_1}$
so
\[
 Q_1(x,x+e_1) + Q_1(x,x-e_1) + Q_1(x,x) = \frac{q_1^{x_1 + 1}}{q_1^{x_1}} + \frac{q_1^{x_1-1}}{q_1^{x_1}}\Indi_{[x_1 > 0]} -  q_1 - q_1^{-1}\Indi_{[x_1 > 0]},
\]
a quantity equal to zero, and the off diagonal entries are clearly positive.

Under the assumption that H($2n-1$) holds we will define a conservative $Q$-matrix 
$\Genwo$ on $\Sswo=  \{ (x,y) \in \Sxw \times \Sxw : x \prec y\}$ in terms of $Q_{2n-1}$ and prove the intertwining relationship
\[
 Q_{2n} \Lambda = \Lambda \Genwo.
\]
where $\Lambda$ is a Markov kernel.  Expanding the intertwining and summing both sides shows that $\sum_{x^\prime} Q_{2n}(x,x^\prime) = 0$, so we conclude that H($2n$) holds as well.  The step from H($2n$) to H($2n+1$) follows a similar argument.

\subsection{Part I: Iterating from an odd row to an even row}\label{s:poissonwallI}

Suppose H($2n-1$) holds and identify $Q_X \equiv Q_{2n-1}$.  Introduce a $Q$-matrix $\Genwo$ on $\Sswo$ 
with off diagonal entries defined by
\[
\Genwo((x,y),(x^\prime,y^\prime)) = \left\{\begin{array}{ll}
Q_X(x,x \pm e_j), &\; (x^\prime,y^\prime) = (x \pm e_j,y)  \\
Q_X(x,x - e_{i+1}), &\; (x^\prime,y^\prime) = (x - e_{i+1},y - e_i),\; x_{i+1} = y_i \\
Q_X(x,x + e_j), &\; (x^\prime,y^\prime) = (x + e_j,y + e_j),\; x_j = y_j \\
q^{\mp 1}_{n}, &\; (x^\prime,y^\prime) = (x,y \pm e_j) \\
0 &\; \mathrm{otherwise}
\end{array}\right.,
\]
for $(x,y),(x^\prime,y^\prime) \in \Sswo$, $1 \leq i < n$, $1 \leq j \leq n$.
The diagonal entry $-\Genwo((x,y),(x,y))$ is given by
 \begin{equation}\label{e:oddrateleaving}
 \sum_{i=1}^{n-1} \left(q_i + q_i^{-1}\right) + q_n + q^{-1}_n \Indi_{[x_1 > 0]}
+ \sum_{i=1}^{n-1}\left( q^{-1}_n\Indi_{[y_i < x_{i+1}]} + q_n \Indi_{[y_i > x_i]} \right) + 
 q_n\Indi_{[y_n > x_n]} + q^{-1}_n,
\end{equation}
so under the assumption that $Q_X$ is conservative, $\Genwo$ is also conservative.

Define $m:\Sswo \to [0,1]$ by
\[
m(x,y) = q_{n}^{|x| - |y|} \frac{Sp^{2n-1}_x(q)}{Sp^{2n}_y(q)}.
\]
Note that the geometric factor is now $q_{n}^{|x| - |y|}$ 
instead of the usual $q_{n}^{|y| - |x|}$.
By definition \eqref{e:oddsympschur},

\begin{eqnarray*}
Sp^{2n}_z(q) & = & \sum_{x \in \SympKone_{2n}(z)} w^q_{2n}(q)  \\
 & = & \sum_{x \in \SympKone_{2n}(z)} q_n^{|x^{2n-1}| - |x^{2n}|} q_n^{|x^{2n-1}| - |x^{2n-2}|}\prod_{i=1}^{n-1} q_i^{2|x^{2i-1}| - |x^{2i}| - |x^{2i-2}|} \\
& = &  \sum_{z^\prime \in \Sxw : z^\prime \prec z} q_n^{|z^\prime| - |z|}\sum_{x \in \SympKone_{2n-1}(z^\prime)} q_n^{|x^{2n-1}| - |x^{2n-2}|}\prod_{i=1}^{n-1} q_i^{2|x^{2i-1}| - |x^{2i}| - |x^{2i-2}|} \\
&= & \sum_{z^\prime \in \Sxw : z^\prime \prec z} q_n^{|z^\prime| - |z|} Sp^{2n-1}_{z^\prime}(q).
\end{eqnarray*}
Hence, $m$ gives a Markov kernel $\Lambda$ from $\Sxw$ to $\Sswo$ defined by
\[
 \Lambda(y,(x^\prime, y^\prime)) = m(x^\prime, y^\prime) \Indi_{[y^\prime = y]}.
\]

We then have

\begin{theorem}\label{t:wallodd}
Assume $Q_{2n-1}$ is a conservative $Q$-matrix and $(X(t),Y(t); t \geq 0)$ is a Markov process with $Q$-matrix $\Genwo$ and initial distribution $\Lambda(y,\cdot)$ for some $y \in \Sxw$.
Then $Q_{2n}$ is a conservative $Q$-matrix and $(Y(t); t \geq 0)$ is distributed as a Markov process with $Q$-matrix $Q_{2n}$, started from $y$.
\end{theorem}

Suppose $Q_Y \equiv Q_{2n}$, then as usual we prove an intertwining relationship
\[
 Q_Y \Lambda = \Lambda \Genwo.
\]

This is equivalent to
\begin{equation}\label{e:QYintertwiningwo}
Q_Y(y,y^\prime) = \sum_{x \prec y} \frac{m(x,y)}{m(x^\prime,y^\prime)} \Genwo((x, y),(x^\prime, y^\prime)), \;\; y \in \Sxw, \; (x^\prime, y^\prime) \in \Sswo.
\end{equation}
where the sum is over $x \in \Sxw$ such that $(x,y) \in \Sswo$.

Particles may take unit steps in either direction so we need to 
check the equality \eqref{e:QYintertwiningwo} holds for
$y = y^\prime$, $y = y^\prime + e_j$ and $y = y^\prime - e_j$ 
for some $1 \leq j \leq n$.

Let us first consider the case $y = y^\prime$.  When $x = x^\prime$,
$-\Genwo((x, y),(x^\prime, y^\prime))$ is the rate of leaving $(x^\prime,y^\prime)$ and is given by
\eqref{e:oddrateleaving}.  The only other possible values of $x$ in the summation for which the summand is non-zero are $x = x^\prime \pm e_i$, $1 \leq i \leq n$.  For such values (i.e.
if $(x^\prime \pm e_i, y) \in \Sswo$), the summand is
\[
\frac{m(x^\prime \pm e_i,y^\prime)}{m(x^\prime,y^\prime)} Q_X(x^\prime \pm e_i, x^\prime)
= \frac{q_n^{|x^\prime \pm e_i|-|y^\prime|}}{q_n^{|x^\prime|-|y^\prime|}} 
\frac{Sp^{2n-1}_{x^\prime \pm e_i}(q)}{ Sp^{2n}_{y^\prime}(q)}
\frac{Sp^{2n}_{y^\prime}(q)}{Sp^{2n-1}_{x^\prime}(q)} 
\frac{Sp^{2n-1}_{x^\prime}(q)}{Sp^{2n-1}_{x^\prime \pm e_i}(q)}, 
\]
which is a rather fancy way of writing $q_n^{\pm 1}$.  But, for $(x^\prime,y^\prime) \in \Sswo$,
\begin{itemize}
 \item $(x^\prime + e_i, y^\prime) \in \Sswo$ only if $x^\prime_i < y^\prime_i$
 \item $(x^\prime - e_i, y^\prime) \in \Sswo$, $i > 1$ only if $x^\prime_i > y^\prime_{i-1}$
 \item $(x^\prime - e_1, y^\prime) \in \Sswo$, only if $x^\prime_1 > 0$.
\end{itemize}
So
\begin{eqnarray*}
\sum_{x \prec y, x \neq x^\prime} \frac{m(x,y)}{m(x^\prime,y^\prime)} \Genwo((x, y),(x^\prime, y^\prime)) & = &\sum_{i=1}^{n-1} \left(q_n \Indi_{[x^\prime_i < y_i]} + q^{-1}_n \Indi_{[x^\prime_{i+1} > y_i]} \right) + \\
 && \qquad q_n \Indi_{[x^\prime_n < y_n]} + q^{-1}_n \Indi_{[x^\prime_1 > 0]}
\end{eqnarray*}

On subtracting the rate of leaving $-\Genwo((x^\prime,y^\prime),(x^\prime,y^\prime))$ defined in \eqref{e:oddrateleaving} 
we find that the indicator functions
all cancel and the right hand side of \eqref{e:QYintertwiningwo} is 
\[
 \sum_{i=1}^{n} \left(q_i + q_i^{-1}\right),
\]
which is equal to the left hand side.

Next we consider the case that $y^\prime = y - e_i \in \Sxw$.  If $i = n$, 
the only possibility is that the $Y$ particle jumped by itself.  When $i < n$, the only 
possibilities are that the $i^{th}$ component of $Y$ was pushed by 
the $(i+1)^{th}$ component of $X$ (i.e. $x = x^\prime + e_{i+1}$) or it jumped by its own volition (i.e. $x = x^\prime$).
The former only occurs if $y^\prime_i = x^\prime_{i+1}$, while the latter 
can only occur if $y^\prime_i < x^\prime_{i+1}$, inducing a natural partition 
on the values we have to check the intertwining on.  When $y^\prime = y - e_i$, $y^\prime_i < x^\prime_{i+1}$, $i < n$, or $i = n$,
the right hand side of \eqref{e:QYintertwiningwo} is 
\[
 \frac{m(x^\prime,y^\prime + e_i)}{m(x^\prime,y^\prime)} \Genwo((x^\prime,y^\prime + e_i),(x^\prime,y^\prime))
= \frac{q_n^{|x^\prime|-|y^\prime + e_i|}}{q_n^{|x^\prime|-|y^\prime|}} 
\frac{Sp^{2n-1}_{x^\prime}(q)}{ Sp^{2n}_{y^\prime + e_i}(q)}
\frac{Sp^{2n}_{y^\prime}(q)}{Sp^{2n-1}_{x^\prime}(q)} q_n= \frac{Sp^{2n}_{y^\prime}(q)}{ Sp^{2n}_{y^\prime + e_i}(q)}.
\]

When $y^\prime = y - e_i$, $x^\prime_{i+1} = y^\prime_i$, $i < n$, the 
sum on the right hand side of the intertwining involves a single term,
\[
  \frac{m(x^\prime + e_{i+1},y^\prime + e_i)}{m(x^\prime,y^\prime)} Q_X(x^\prime + e_{i+1},x^\prime).
\]
Using the definitions of $m$ and $Q_X$ shows this summand is
\[
\frac{q_n^{|x^\prime + e_{i+1}|-|y^\prime + e_i|}}{q_n^{|x^\prime|-|y^\prime|}} 
\frac{Sp^{2n-1}_{x^\prime+e_{i+1}}(q)}{ Sp^{2n}_{y^\prime + e_i}(q)}
\frac{Sp^{2n}_{y^\prime}(q)}{Sp^{2n-1}_{x^\prime}(q)} 
\frac{Sp^{2n-1}_{x^\prime}(q)}{Sp^{2n-1}_{x^\prime + e_{i+1}}(q)}.
\]
Both of these quantities are equal to $Q_Y(y^\prime+e_i,y)$.

Finally we consider the case $y^\prime = y + e_i$, $1 \leq i \leq n$.  As in the 
previous case, the dichotomy $x^\prime_i = y^\prime_i$ and $x^\prime_i < y^\prime_i$ 
divides the possible values of $x$ in the summation into two cases, each of which 
having only one term contributing to the sum.
When $x^\prime_i = y^\prime_i$, the $i^{th}$ $Y$ particle must have been pushed, and
\[
  \frac{m(x^\prime - e_i,y^\prime - e_i)}{m(x^\prime,y^\prime)} Q_X(x^\prime - e_{i},x^\prime)
= \frac{q_n^{|x^\prime - e_{i}|-|y^\prime - e_i|}}{q_n^{|x^\prime|-|y^\prime|}} 
\frac{Sp^{2n-1}_{x^\prime-e_{i}}(q)}{ Sp^{2n}_{y^\prime - e_i}(q)}
\frac{Sp^{2n}_{y^\prime}(q)}{Sp^{2n-1}_{x^\prime}(q)} 
\frac{Sp^{2n-1}_{x^\prime}(q)}{Sp^{2n-1}_{x^\prime - e_{i}}(q)}.
\]
Simplifying the expression on the right hand side by cancelling common factors
in the numerator and denominator reveal it to be simply $Q_Y(y^\prime-e_i,y)$.

On the other hand, when $x^\prime_i < y^\prime_i$ the $i^{th}$ $Y$ particle
cannot have been pushed so the right hand side of the intertwining 
\eqref{e:QYintertwiningwo} is
\[
  \frac{m(x^\prime,y^\prime - e_i)}{m(x^\prime,y^\prime)} \Genwo((x^\prime,y^\prime - e_i),(x^\prime,y^\prime))
= \frac{q_n^{|x^\prime|-|y^\prime - e_i|}}{q_n^{|x^\prime|-|y^\prime|}} 
\frac{Sp^{2n-1}_{x^\prime}(q)}{ Sp^{2n}_{y^\prime - e_i}(q)}
\frac{Sp^{2n}_{y^\prime}(q)}{Sp^{2n-1}_{x^\prime}(q)} q^{-1}_n= \frac{Sp^{2n}_{y^\prime}(q)}{ Sp^{2n}_{y^\prime - e_i}(q)}.
\]
The proof of the intertwining relationship is concluded by noting that this
is $Q_Y(y^\prime - e_i, y)$ as required.

Now, summing both sides of the intertwining
\[
 \sum_{(x^\prime,y^\prime)} Q_Y(y,y^\prime)m(x^\prime,y^\prime) = \sum_{(x^\prime,y^\prime)} \sum_{x \prec y} m(x,y) \Genwo((x, y),(x^\prime, y^\prime))
\]
over all pairs in $(x^\prime, y^\prime)$ in $\Sswo$ shows that $Q_Y$ is conservative
as $\sum_{x^\prime} m(x^\prime,y^\prime) = 1$ and $\sum_{(x^\prime,y^\prime)} \Genwo((x, y),(x^\prime, y^\prime)) = 0$.

We then apply lemma \ref{l:qmatintimpkernint} to recover the rest of the theorem.

\subsection{Part II: Iterating from an even row to an odd}\label{s:poissonwallII}

Suppose $Q_X \equiv Q_{2n}$, $Q_Y \equiv Q_{2n+1}$ and $\Genw$ is a conservative $Q$-matrix $\Genw$ on $\Ssw = \{ (x,y) \in \Sxw \times \Syw : x \preceq y\}$ with off diagonal entries given by
\[
\Genw((x,y),(x^\prime,y^\prime)) = \left\{\begin{array}{ll}
Q_X(x,x \pm e_i), &\; (x^\prime,y^\prime) = (x \pm e_i,y)  \\
Q_X(x,x - e_i), &\; (x^\prime,y^\prime) = (x - e_i,y - e_i),\; x_i = y_i \\
Q_X(x,x + e_i), &\; (x^\prime,y^\prime) = (x + e_i,y + e_{i+1}),\; x_i = y_{i+1} \\
q^{\pm 1}_{n+1}, &\; (x^\prime,y^\prime) = (x,y \pm e_j) \\
0 &\; \mathrm{otherwise}
\end{array}\right.,
\]
for $(x,y),(x^\prime,y^\prime) \in \Ssw$, $1 \leq i \leq n$, $1 \leq j \leq n+1$.
The diagonal $-\Genw((x, y),(x, y))$ is given by
\begin{equation}\label{e:evenrateleaving}
 \sum_{i=1}^n (q_i + q_i^{-1}) + \sum_{i=1}^{n}(q_{n+1}\Indi_{[y_i < x_i]} + q^{-1}_{n+1}\Indi_{[y_{i+1} > x_i]} ) + q_{n+1} + q^{-1}_{n+1} \Indi_{[y_1 > 0]}.
\end{equation}
Hence $\Genw$ is conservative if $Q_X$ is.

From definition \eqref{e:oddsympschur} we
calculate
\begin{eqnarray*}
Sp^{2n+1}_z(\tilde q) & = & \sum_{x \in \SympKone_{2n+1}(z)} w^q_{2n+1}(\tilde q)  \\
 & = & \sum_{x \in \SympKone_{2n+1}(z)} q_{n+1}^{|x^{2n+1}| - |x^{2n}|} \prod_{i=1}^{n} q_i^{2|x^{2i-1}| - |x^{2i}| - |x^{2i-2}|} \\
 & = & \sum_{z^\prime \in \Sxw : z^\prime \preceq z} q_{n+1}^{|z| - |z^\prime|} \sum_{x \in \SympKone_{2n}(z^\prime)}  \prod_{i=1}^{n} q_i^{2|x^{2i-1}| - |x^{2i}| - |x^{2i-2}|} \\
&= & \sum_{z^\prime \in \Sxw : z^\prime \preceq z} q_{n+1}^{|z| - |z^\prime|} Sp^{2n}_{z^\prime}(q).
\end{eqnarray*}

So the function $m:\Ssw \to [0,1]$ given by
\[
m(x,y) = q_{n+1}^{|y| - |x|} \frac{Sp^{2n}_x(q)}{Sp^{2n+1}_y(\tilde q)}.
\]
induces a Markov kernel from $\Syw$ to $\Ssw$,
\[
 \Lambda(y,(x^\prime, y^\prime)) = m(x^\prime, y^\prime) \Indi_{[y^\prime = y]}.
\]

Our final theorem is

\begin{theorem}\label{t:walleven}
Assume $Q_{2n}$ is a conservative $Q$-matrix and suppose $(X(t),Y(t); t \geq 0)$ is a Markov process with $Q$-matrix 
$\Genw$ and initial distribution $\Lambda(y,\cdot)$ for some $y \in \Syw$.
Then $Q_{2n+1}$ is a conservative $Q$-matrix and $(Y(t); t \geq 0)$ is 
distributed as a Markov process with $Q$-matrix $Q_{2n+1}$, started from $y$.
\end{theorem}

The intertwining via $\Lambda$ is equivalent to
\begin{equation}\label{e:QYintertwiningw}
Q_Y(y,y^\prime) = \sum_{x \preceq y} \frac{m(x,y)}{m(x^\prime,y^\prime)} \Genw((x, y),(x^\prime, y^\prime)), \;\; y \in \Syw, \; (x^\prime, y^\prime) \in \Ssw,
\end{equation}
where we sum over $x$ such that $(x,y) \in \Ssw$.

We only need to check \eqref{e:QYintertwiningw} holds for $y$ of the form
$y = y^\prime$, $y = y^\prime \pm e_j$ for $1 \leq j \leq n + 1$ as both sides
vanish otherwise.

Again we start with the case $y^\prime = y$.  When $x = x^\prime$, the rate of leaving $-\Genw((x^\prime, y^\prime),(x^\prime, y^\prime))$ is given by 
\eqref{e:evenrateleaving}.  The only other possible values of $x$ 
for which the summand is non-zero are $x = x^\prime \pm e_i$ for 
$1 \leq i \leq n$.  For such $x$ values satisfying $(x,y) \in \Ssw$, the definitions of $m$ and $Q_X$ give
\[
 \frac{m(x^\prime \pm e_i,y^\prime)}{m(x^\prime,y^\prime)} Q_X(x^\prime \pm e_i,x^\prime) = 
\frac{q_{n+1}^{|y^\prime| - |x^\prime \pm e_i|}}{q_{n+1}^{|y^\prime| - |x^\prime|}}
\frac{Sp^{2n}_{x^\prime \pm e_i}(q)}{Sp^{2n+1}_{y^\prime}(\tilde q)}
\frac{Sp^{2n+1}_{y^\prime}(\tilde q)}{Sp^{2n}_{x^\prime}(q)}
\frac{Sp^{2n}_{x^\prime}(q)}{Sp^{2n}_{x^\prime \pm e_i}(q)},
\]
which is equal to $q_{n+1}^{\mp 1}$.  But, for $(x^\prime,y^\prime) \in \Ssw$,
\begin{itemize}
 \item $(x^\prime + e_i, y^\prime) \in \Ssw$ only if $x^\prime_i < y^\prime_{i+1}$ and
 \item $(x^\prime - e_i, y^\prime) \in \Ssw$ only if $x^\prime_i > y^\prime_{i}$.
\end{itemize}

So,
\[
\sum_{x \prec y, x \neq x^\prime} \frac{m(x,y)}{m(x^\prime,y^\prime)} \Genwo((x, y),(x^\prime, y^\prime)) = \sum_{i=1}^{n} \left(q^{-1}_{n+1} \Indi_{[x^\prime_i < y_{i+1}]} + q_{n+1} \Indi_{[x^\prime_{i} > y_i]} \right).
\]

If we now subtract the rate of leaving \eqref{e:evenrateleaving} we
find that  at $y = y^\prime$ the right hand side of \eqref{e:QYintertwiningw} is equal to
\[
 - \left(\sum_{i=1}^n (q_i + q_i^{-1}) + q_{n+1} + q^{-1}_{n+1} \Indi_{[y^\prime_1 > 0]}\right),
\]
which is equal to $Q_Y(y^\prime,y^\prime)$.

The remaining cases are $y = y^\prime \pm e_i$ for some $1 \leq i \leq n+1$.  Let
us deal with $y^\prime = y - e_i$.  If $i = n+1$, this case corresponds to a leftward jump in the rightmost $Y$ particle, a situation that cannot arise 
through pushing by an $X$ particle. 
If $i < n+1$, then the jump arose by pushing if $x^\prime_i = y^\prime_i$,
while if $x^\prime_i > y^\prime_i$ then the $Y$ particle jumped 
by its own volition.  In the case of pushing ($i < n+1$, $x^\prime_i = y^\prime_i$),
familiar calculations show
\[
 \frac{m(x^\prime + e_i, y^\prime + e_i)}{m(x^\prime, y^\prime)}Q_X(x^\prime+ e_i,x^\prime)  = \frac{Sp^{2n+1}_{y^\prime}(\tilde q)}{Sp^{2n+1}_{y^\prime + e_{i}}(\tilde q)} = Q_Y(y,y^\prime).
\]
In the case of no pushing, i.e. $i < n+1$ and $x^\prime_i > y^\prime_i$ or $i = n+1$,
the summand is
\[
 \frac{m(x^\prime, y^\prime + e_i)}{m(x^\prime, y^\prime)} \Genw((x^\prime, y^\prime + e_i),(x^\prime, y^\prime)) = \frac{Sp^{2n+1}_{y^\prime}(\tilde q)}{Sp^{2n+1}_{y^\prime + e_{i}}(\tilde q)} = Q_Y(y,y^\prime).
\]

Finally we consider the case $y^\prime = y + e_i$, $1 \leq i \leq n+1$ corresponding
to a rightward jump in the $i^{th}$ $Y$ particle.
For $i > 1$, consider the dichotomy $x^\prime_{i-1} = y^\prime_{i}$ 
or $x^\prime_{i-1} < y^\prime_{i}$, corresponding to 
the $i^{th}$ $Y$ particle being pushed upwards by the $(i-1)^{th}$ $X$ particle and a free jump respectively.  The case $i = 1$
corresponds to the leftmost $Y$ particle jumping rightwards,
an event that cannot arise as a result of pushing.
In the case of pushing, i.e. $i >1$ and $x^\prime_{i-1} = y^\prime_{i}$, the summand
is equal to
\[
\frac{m(x^\prime - e_{i-1}, y^\prime - e_i)}{m(x^\prime, y^\prime)} \Genw((x^\prime - e_{i-1}, y^\prime - e_i),(x^\prime, y^\prime)).
\]
Using the definitions of $\Genw$ and $m$, this is
\[
\frac{q_{n+1}^{|y^\prime - e_i| - |x^\prime - e_{i-1}|}}{q_{n+1}^{|y^\prime| - |x^\prime|}}\frac{Sp^{2n}_{x^\prime - e_{i-1}}(q)}{Sp^{2n+1}_{y^\prime - e_{i}}(\tilde q)}\frac{Sp^{2n+1}_{y^\prime}(\tilde q)}{Sp^{2n}_{x^\prime}(q)}Q_X(x^\prime - e_{i-1},x^\prime) = \frac{Sp^{2n+1}_{y^\prime}(\tilde q)}{Sp^{2n+1}_{y^\prime - e_{i}}(\tilde q)} = Q_Y(y,y^\prime),
\]
as required.

If $y_i^\prime > x^\prime_{i-1}$ ($i > 1$) or $i = 1$, then the $i^{th}$ $Y$ particle jumped
of its own accord and the only term in the summation is
\[
\frac{m(x^\prime, y^\prime - e_i)}{m(x^\prime, y^\prime)} \Genw((x^\prime, y^\prime - e_i),(x^\prime, y^\prime)) = \frac{Sp^{2n+1}_{y^\prime}(\tilde q)}{Sp^{2n+1}_{y^\prime - e_{i}}(\tilde q)} = Q_Y(y,y^\prime).
\]

This concludes the verification of the intertwining relationship and the theorem follows.

\appendix

\section{A lemma on intertwinings of $Q$-matrices}

\begin{lemma}\label{l:qmatintimpkernint}
Suppose that $L$ and $L^\prime$ are uniformly bounded conservative $Q$-matrices 
on discrete spaces $U$ and $V$ that are intertwined by 
a Markov kernel $\Lambda:U \times V \to [0,1]$ from $U$ to $V$, i.e.
\[
 L \Lambda = \Lambda L^\prime.
\]
Then the transition kernels for the Markov processes with $Q$-matrices $L$ and $L^\prime$ are also intertwined.
\end{lemma}

Note that we always use the lemma with $U = W^\prime, V = W \times W^\prime$ 
where $W \subset \Z^n$ and either $W^\prime \subset \Z^n$ or $W^\prime \subset \Z^{n+1}$. 
Our Markov kernel $\Lambda$ is always such that $\Lambda(u,(v,v^\prime)) > 0$ only 
if $u = v^\prime$, but this, of course, is not necessary for the lemma.

\begin{proof}
The intertwining relationship $L \Lambda = \Lambda L^\prime$ may be written
\[
\sum_{\tilde v \in V} \Lambda(u, \tilde v) L^\prime(\tilde v, v) = \sum_{\tilde u \in U} L(u, \tilde u) \Lambda(\tilde u, v) , \quad u \in U, v \in V.
\]

Let $(p_t;t \geq 0)$ denote the transition kernels for the Markov process corresponding to $Q$-matrix $L$
and fix $u_0 \in U$.  Multiplying both sides of the expanded intertwining relationship 
above by $p_t(u_0,u)$ and summing over $u \in U$ gives
\begin{equation}\label{e:expandedintertwining}
 \sum_{u \in U }p_t(u_0,u) \sum_{\tilde v \in V} \Lambda(u, \tilde v) L^\prime(\tilde v, v) = \sum_{u \in U }p_t(u_0,u) \sum_{\tilde u \in U} L(u, \tilde u) \Lambda(\tilde u, v) , \quad v \in V.
\end{equation}
Now, let $c \in \R$ be a uniform bound for the absolute values of the entries of $L$ and $L^\prime$.
Then 
\[
|p_t(u_0,u) \Lambda(u, \tilde v) L^\prime(\tilde v, v)| \leq c \Lambda(u, \tilde v) p_t(u_0,u), 
\]
so the double sum on the left hand side is absolutely convergent.
Also, 
\[
 \sum_{\tilde u \in U}|p_t(u_0,u) \Lambda(\tilde u, v)  L(u, \tilde u)| \leq \sum_{\tilde u \in U}|p_t(u_0,u) L(u, \tilde u)| \leq 2c p_t(u_0,u)
\]
and the same conclusion holds for the double sum on the right hand side.
So, we may exchange the order of the sums on both sides to give
\[
 \sum_{\tilde v \in V} L^\prime(\tilde v, v) \sum_{u \in U }p_t(u_0,u)  \Lambda(u, \tilde v)  = \sum_{\tilde u \in U} \Lambda(\tilde u, v) \sum_{u \in U } p_t(u_0,u)  L(u, \tilde u)  , \quad v \in V.
\]

Now, as $(p_t; t \geq 0)$ is the transition kernel corresponding to the Markov process with $Q$ matrix $L$, 
it satisfies the Kolmogorov forward equation
\[
 \frac{d}{dt} p_t(u_0, \tilde u) = \sum_{u \in U } p_t(u_0,u)  L(u, \tilde u).
\]
Let us define $q_t(v) = \sum_{\tilde u \in U } p_t(u_0, \tilde u)  \Lambda(\tilde u, v)$ for $v \in V$.  

We may differentiate the summation term by term in $t$ using 
Fubini's theorem and the absolute bounds on the
summands discussed above.  Hence, the right hand side of \eqref{e:expandedintertwining} is simply $\frac{d}{dt} q_t(v)$.

Then, using the definition of $q_t$ in the left hand side of \eqref{e:expandedintertwining}, we see that 
\begin{equation}\label{e:primeforwardeq}
\frac{d}{dt} q_t(v) = \sum_{\tilde v \in V} L^\prime(\tilde v, v) q_t(\tilde v), \quad t \geq 0.
\end{equation}

Now let $(p^\prime_t; t \geq 0)$ denote the transition kernels of the Markov process with $Q$-matrix $L^\prime$, and 
\[
 p^\prime_t(v) = \sum_{\tilde v} q_0(\tilde v) p^\prime_t(\tilde v, v) =  \sum_{\tilde v \in V} \Lambda(u_0, \tilde v) p^\prime_t(\tilde v, v).
\]
Then $p^\prime_0(v) = q_0(v)$ for all $v \in V$ and $p^\prime_t$ also satisfies the forward equation \eqref{e:primeforwardeq} in $L^\prime$.

But when the rates are uniformly bounded there is exactly one solution to the forward differential equation with the same boundary conditions as $q_t$ so $q_t(v) = p^\prime_t(v)$ for all $t \geq 0$ and $v \in V$.

By definition of $q_t(v)$, we then have
\[
\sum_{\tilde u \in U } p_t(u_0, \tilde u)  \Lambda(\tilde u, v) = \sum_{\tilde v \in V} \Lambda(u_0, \tilde v) p^\prime_t(\tilde v, v),
\]
and since the argument holds for arbitrary $u_0 \in U$ we're done.

\end{proof}


\end{document}